\documentclass[12pt]{amsart}

\usepackage{epsfig}
\usepackage{amssymb}
\usepackage{enumerate}
\usepackage{subfigure}
\usepackage{setspace}
\usepackage{graphicx}
\usepackage{epstopdf}

\newcommand{\A}{\mathcal{A}}

\newcommand{\RR}{\mathbb{R}}
\newcommand{\ZZ}{\mathbb{Z}}
\newcommand{\FF}{{\mathbb F}}
\newcommand{\Cox}{{\sf Cox}}
\newcommand{\Shi}{{\sf Shi}}
\newcommand{\Ish}{{\sf Ish}}
\newcommand{\Rec}{{\sf Rec}}
\renewcommand{\L}{{\mathfrak L}}

\newcommand{\Stir}{{\sf Stir}}
\renewcommand{\a}{{\bf a}}
\renewcommand{\b}{{\bf b}}

\numberwithin{equation}{section}

\newtheorem{thm}{Theorem}[section]

\newtheorem{cor}[thm]{Corollary}
\newtheorem{lem}[thm]{Lemma}

\newtheorem*{unmainthm}{Main Theorem}
\theoremstyle{definition}
\newtheorem{defn}{Definition}[section]

\numberwithin{figure}{section}

\begin{document}

\title{The Shi Arrangement and the Ish Arrangement}

\author{Drew Armstrong}
\address{Drew Armstrong, Department of Mathematics, University of Miami, Coral
Gables, FL, 33146}
\email{armstrong@math.miami.edu}
\author{Brendon Rhoades}
\address{Brendon Rhoades, Department of Mathematics, MIT, Cambridge, MA, 02139}
\email{brhoades@math.mit.edu}



\bibliographystyle{../dart}

\date{\today}

\maketitle

\begin{abstract}
This paper is about two arrangements of hyperplanes. The first --- the {\sf Shi arrangement} --- was introduced by Jian-Yi Shi \cite[Chapter 7]{Shi} to describe the Kazhdan-Lusztig cells in the affine Weyl group of type $A$. The second --- the {\sf Ish arrangement} --- was recently defined by the first author \cite{Armqt} who used the two arrangements together to give a new interpretation of the $q,t$-Catalan numbers of Garsia and Haiman. In the present paper we will define a mysterious ``combinatorial symmetry'' between the two arrangements and show that this symmetry preserves a great deal of information. For example, the Shi and Ish arrangements share the same characteristic polynomial, the same numbers of regions, bounded regions, dominant regions, regions with $c$ ``ceilings'' and $d$ ``degrees of freedom'', etc. Moreover, all of these results hold in the greater generality of ``deleted'' Shi and Ish arrangements corresponding to an arbitrary subgraph of the complete graph. Our proofs are based on nice combinatorial labelings of Shi and Ish regions and a new set
partition-valued statistic on these regions.
\end{abstract}

\section{Introduction}
A {\sf hyperplane arrangement} is a finite collection of affine hyperplanes in Euclidean space. Some of the nicest arrangements come from the reflecting hyperplanes of Coxeter groups. In particular, the {\sf Coxeter arrangement of type $A$} (also known as the {\sf braid arrangement}) is 
the arrangement in $\RR^n$
defined by
\begin{equation*}
\Cox(n):=\{ x_i-x_j = 0 : 1\leq i<j\leq n \}.
\end{equation*}
Here $\{x_1,\ldots,x_n\}$ are the standard coordinate functions on $\RR^n$.

\begin{figure}
\centering
\includegraphics[scale=.6]{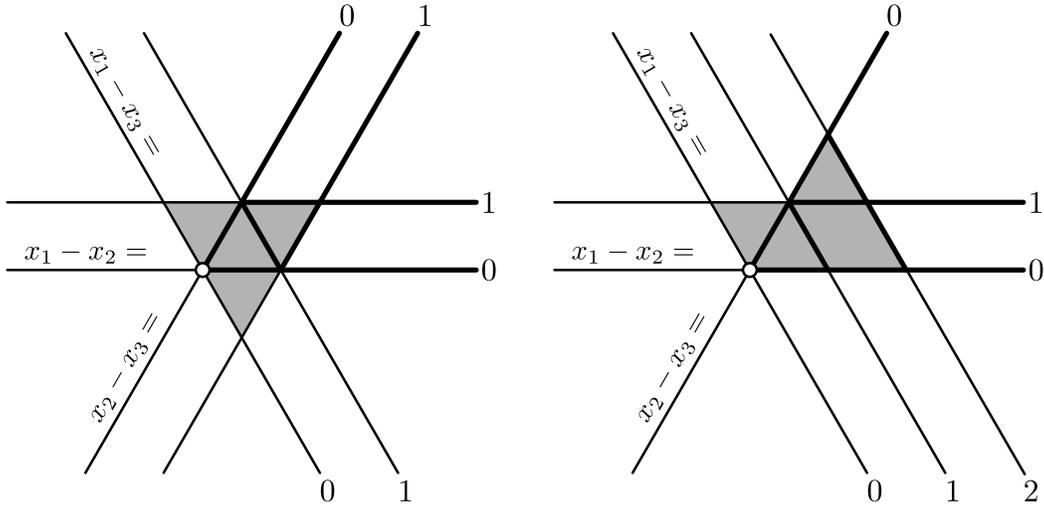}
\caption{The arrangements $\Shi(3)$ (left) and $\Ish(3)$ (right)}
\label{fig:shiishthree}
\end{figure}

Postnikov and Stanley \cite{PostnikovStanley} introduced the idea of a {\sf deformation} of the Coxeter arrangement --- this is an affine arrangement each of whose hyperplanes is parallel to some hyperplane of the Coxeter arrangement.  
In the present paper we will study two specific deformations of the Coxeter arrangement and we will observe a deep similarity between them. The first is the {\sf Shi arrangement} which was one of Postnikov and Stanley's motivating examples:
\begin{equation*}
\Shi(n):=\Cox(n) \cup \{ x_i-x_j=1: 1\leq i<j\leq n\}.
\end{equation*}
This arrangement was defined by Jian-Yi Shi \cite[Chapter 7]{Shi} in
the study of the Kazhdan-Lusztig cellular structure of
the affine Weyl group of type A. The second is the {\sf Ish arrangement}, recently defined by the first author \cite{Armqt}:
\begin{equation*}
\Ish(n):= \Cox(n) \cup \{ x_1-x_j = i : 1\leq i<j\leq n\}.
\end{equation*}
He used the Shi and Ish arrangements to give a new description of the $q,t$-Catalan numbers of Garsia and Haiman in terms of the affine Weyl group of type $A$. Figure \ref{fig:shiishthree} displays the arrangements $\Shi(3)$ and $\Ish(3)$. (Note that the normals to the hyperplanes of either $\Shi(n)$ or $\Ish(n)$ span the hyperplane $x_1+x_2+\cdots +x_n=0$. Hence we will always draw their restrictions to this space.)

\begin{figure}
\centering
\includegraphics[scale=.6]{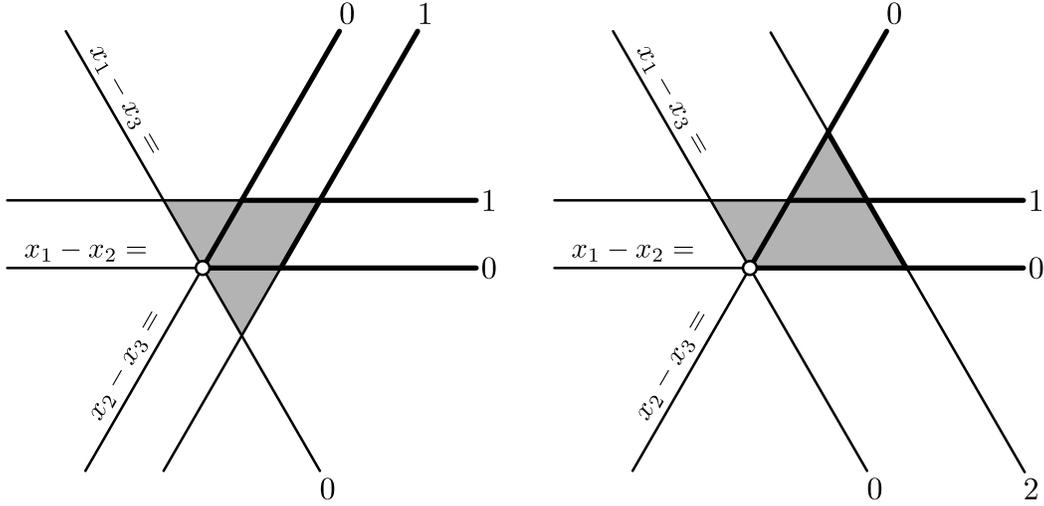}
\caption{The arrangements $\Shi(G)$ (left) and $\Ish(G)$ (right) corresponding to the ``chain'' $G=\{12,23\}\subseteq\binom{[3]}{2}$}
\label{fig:shiishdelete}
\end{figure}

The heart of this paper is the following correspondence between Shi and Ish hyperplanes. The correspondence  is natural to state but we find it geometrically mysterious. We will call this a ``combinatorial symmetry'':
\begin{equation*}
\boxed{ x_i-x_j=1 \longleftrightarrow x_1-x_j = i \quad\text{ for }\quad 1\leq i<j\leq n}
\end{equation*}
This symmetry allows us to define {\sf deleted} versions of the Shi and Ish arrangements. Let $\binom{[n]}{2}$ denote the set of pairs $ij$ satisfying $1\leq i<j\leq n$ and consider a simple loopless graph $G\subseteq\binom{[n]}{2}$. The {\sf deleted Shi and Ish arrangement} are defined as follows:
\begin{align*}
\Shi(G) &:= \Cox(n) \cup \{ x_i-x_j=1 : ij \in G\},\\
\Ish(G) &:= \Cox(n) \cup \{x_1-x_j=i : ij\in G\}.
\end{align*}
The arrangement $\Shi(G)$ was first considered by Athanasiadis \cite{Athanasiadis}. Note that $\Shi(G)$ (resp. $\Ish(G)$) interpolates between the Coxeter arrangement and the Shi (resp. Ish) arrangement. That is, if $\emptyset\in\binom{[n]}{2}$ is the ``empty'' graph and $K_n=\binom{[n]}{2}$ is the ``complete'' graph, we have
\begin{equation*}
\Shi(\emptyset)=\Ish(\emptyset)=\Cox(n), \quad \Shi(K_n)=\Shi(n)\quad \text{and}\quad \Ish(K_n)=\Ish(n).
\end{equation*}
Figure \ref{fig:shiishdelete} displays the arrangements $\Shi(G)$ and $\Ish(G)$ corresponding to the ``chain'' $G=\{12,23\}\subseteq\binom{[3]}{2}$.

In order to state our main results right away we need a few definitions.

Let $\A$ be either $\Shi(G)$ or $\Ish(G)$. The connected components of $\RR^n-\cup_{H\in \A} H$ are called {\sf regions}. We say that a region is {\sf dominant} if it lies within the {\sf dominant cone}, defined by the coordinate inequalities
\begin{equation*}
x_1>x_2>\cdots >x_n.
\end{equation*}
The topological closure $\bar{R}$ of a region $R$ is decomposed by the arrangement $\A$ into {\sf faces} of various dimensions. We say that the hyperplane $H$ is a {\sf wall} of $R$ if it is the affine span of a codimension-$1$ face of $R$. The wall $H$ is called a {\sf ceiling} if 
$H$ does not contain the origin and if
the region $R$ and the origin lie in the same half-space of $H$. Since every region $R$ is convex, it determines a {\sf recession cone} (which is closed under non-negative linear combinations):
\begin{equation*}
\Rec(R):=\{ v\in\RR^n : v + R \subseteq R \}.
\end{equation*}
Note that the region $R$ is bounded if and only if
$\Rec(R) = 0$. We call the dimension of $\Rec(R)$ the number of {\sf degrees of freedom} of $R$.

Finally, let $\L(\A)$ denote the collection of intersections of hyperplanes from $\A$, partially ordered by {\em reverse}-inclusion of subspaces:
\begin{equation*}
\L(\A):=\left\{ \cap_{H\in S} H : S\subseteq \A\right\}.
\end{equation*}
This poset has the structure of a geometric semilattice (see \cite{WachsWalker}) with a unique minimum element $\RR^n$ (corresponding to the empty intersection). The {\sf characteristic polynomial} (or {\sf chromatic polynomial}) $\chi_\A(p)\in\ZZ[p]$ of the arrangement $\A$ is defined by
\begin{equation*}
\chi_{\A}(p) = \sum_{X \in \mathfrak{L}(\A)} \mu(\RR^n,X) p^{\mathrm{dim}(X)},
\end{equation*}
where $\mu: \mathfrak{L}(\A) \times \mathfrak{L}(\A) \rightarrow \ZZ$ is the 
M\"obius function of the poset $\mathfrak{L}(\A)$ (see \cite{StanEC1}).


\begin{unmainthm} Let $G\subseteq\binom{[n]}{2}$ be a simple loopless graph on $n$ vertices; let $c$
and $d$ be nonnegative integers. The deleted Shi and Ish arrangements $\Shi(G)$ and $\Ish(G)$ share the following properties in common:
\begin{enumerate}
\item the characteristic polynomial;
\item the number of {\bf dominant} regions with $c$ ceilings;
\item the number of regions with $c$ ceilings and $d$ degrees of freedom.
\end{enumerate}
\end{unmainthm}
\begin{proof}
Parts (1), (2), (3) are Theorems \ref{th:charpoly}, \ref{thm:dominant}, and \ref{thm:candd}, respectively.
\end{proof}

For example, here are the joint distributions of ceilings ($c$) and degrees of freedom ($d$) for the arrangements in Figures \ref{fig:shiishthree} and \ref{fig:shiishdelete}, respectively.
\begin{center}
\begin{tabular}{cc}
\begin{tabular}{rc}
& \,\,\,\,\,\,\,\,$d$ \\
\begin{tabular}{c} \\ $c$\!\!\!\!\!\!\!\!\!\! \end{tabular} & 
\begin{tabular}{r|ccc}
 & 1 & 2 & 3\\
\hline
0 & & & 6\\
1 & 3 & 6 &\\
2 & 1 & & \\
\end{tabular}\\
\end{tabular}
&
\quad\quad\quad\quad\begin{tabular}{rc}
& \,\,\,\,\,\,\,\,$d$ \\
\begin{tabular}{c} \\ $c$\!\!\!\!\!\!\!\!\!\! \end{tabular} & 
\begin{tabular}{r|ccc}
 & 1 & 2 & 3\\
\hline
0 & & & 6\\
1 & 2 & 4 & \\
2 & 1 & & \\
\end{tabular}\\
\end{tabular}
\end{tabular}
\end{center}

We find it surprising that the symmetry $x_i-x_j=1\leftrightarrow x_1-x_j= i$ preserves so much information. However, there are important properties that it does not preserve. For example, one may observe from Figures \ref{fig:shiishthree} and \ref{fig:shiishdelete} that the intersection poset is not preserved. One can also show that the Tutte polynomials of $\Shi(3)$ and $\Ish(3)$ differ, and that the Orlik-Solomon algebras of $\Shi(G)$ and $\Ish(G)$ are not graded-isomorphic for $G=\{12,23\}$
(even though the equality of characteristic polynomials implies that these algebras have the same Hilbert series). Is there a unifying concept that could simplify the statement of the Main Theorem?

The paper is structured as follows.

In {\bf Section 2} we establish some language for set partitions. We define {\sf $G$-partitions} --- which for the complete graph are just partitions of $[n]=\{1,2,\ldots,n\}$ --- and discuss various kinds: connected, nonnesting. We define the $(\a,\b)$ {\sf endpoint notation} for partitions which seems to be the correct language for comparing Shi and Ish arrangements.

In {\bf Section 3} we show that $\Shi(G)$ and $\Ish(G)$ have the same characteristic polynomial, which has a formula involving $G$-Stirling numbers. This proves part (1) of the Main Theorem. Our tool is the finite field method Crapo and Rota \cite{CrapoRota}. By a standard result of Zaslavsky this implies that $\Shi(G)$ and $\Ish(G)$ share the same numbers of total regions and relatively bounded regions (regions with one degree of freedom).

In {\bf Section 4} we modify a labeling of the regions of $\Shi(G)$ due to Athanasiadis and Linusson \cite{ALShi} and we call the result {\sf Shi ceiling diagrams}. Similarly, we define {\sf Ish ceiling diagrams} for the regions of $\Ish(G)$. We give a bijective proof of part (2) of the Main Theorem by observing that {\bf dominant} regions of $\Shi(G)$ and $\Ish(G)$ correspond to order ideals in isomorphic posets.

In {\bf Section 5} we define the {\sf ceiling partition} for a region of $\Shi(G)$ or $\Ish(G)$. This is a (possibly nesting) $G$-partition that encodes the ceilings of the region. Given a $G$-partition $\pi$ with $k$ blocks and an integer $1\leq d\leq k$, we show that the number of regions of either $\Shi(G)$ or $\Ish(G)$ with ceiling partition $\pi$ and $d$ degrees of freedom is equal to
\begin{equation*}
\frac{d(n-d-1)!(k-1)!}{(n-k-1)!(k-d)!},
\end{equation*}
which proves part (3) of the Main Theorem. This formula is remarkable, and it is new even for the Shi arrangement. The proof of the formula for Ish regions is direct, whereas the proof for Shi regions uses a new formula due to the second author (see \cite{RCat} or Lemma 2.3) which counts nonnesting partitions with a fixed number of connected components and fixed block size multiplicities. This suggests an open problem: Find a bijection between Shi regions and Ish regions with ceiling partition $\pi$ and $d$ degrees of freedom. This bijection {\bf cannot} preserve the property of being dominant, since $\Shi(G)$ and $\Ish(G)$ do {\bf not} share the same number of dominant regions with $d$ degrees of freedom.

We end with an observation:
\begin{quote}
The Ish arrangement is something of a ``toy model'' for the Shi arrangement (and other Catalan objects). That is, for any property $P$ that $\Shi(G)$ and $\Ish(G)$ share, the proof that $\Ish(G)$ satisfies $P$ is easier than the proof that $\Shi(G)$ satisfies $P$.
\end{quote}

\section{Set Partitions}
\label{sec:partitions}
All of the formulas in this paper are phrased in terms of set partitions. In this section we will give some background on these and establish notation. In particular, for each graph $G\subseteq\binom{[n]}{2}$ we will define {\sf $G$-partitions} of the set $[n]=\{1,2,\ldots,n\}$. In the case of the complete graph this corresponds to the usual notion of partitions.

\subsection{The endpoint notation}
\begin{figure}
\centering
\includegraphics[scale=1]{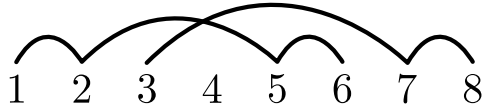}
\caption{A partition of $[8]$ with type $(1,0,1,1,0,0,0,0)$}
\label{fig:arcdiagram}
\end{figure}

We say that $\pi=\{B_1,B_2,\ldots,B_k\}$ is a {\sf partition of $[n]$ into $k$ blocks} if the following disjoint union holds:
\begin{equation*}
[n]=B_1\sqcup B_2\sqcup\cdots\sqcup B_k.
\end{equation*}
The {\sf type} of the partition $\pi$ is the sequence $(r_1,r_2,\ldots,r_n)$ where $r_i$ is the number of blocks of $\pi$ with size $i$. We draw the {\sf arc diagram} of $\pi$ as follows: Place the numbers $1,2,\ldots,n$ on a line and draw an arc between each pair $i<j$ such that
\begin{itemize}
\item $i$ and $j$ are in the same block of $\pi$; and
\item there is no $i<\ell<j$ such that $i,\ell,j$ are in the same block of $\pi$.
\end{itemize}
Figure \ref{fig:arcdiagram} displays the arc diagram for the partition $\{\{1,2,5,6\},\{3,7,8\},\{4\}\}$, which has type $(1,0,1,1,0,0,0,0)$.

In this paper we will use a special notation for partitions, based on the arc diagram. First note that a partition $\pi$ has $n-k$ blocks if and only if its diagram has $k$ arcs. This is because each new arc reduces the number of blocks by one. Now suppose that the arcs of $\pi$ are $a_1b_1, a_2b_2,\ldots, a_kb_k$, with the left endpoints in increasing order: $a_1<a_2<\cdots<a_k$. We will associate $\pi$ with its pair $(\a,\b)$ of {\sf endpoint vectors}:
\begin{equation*}
\a=a_1a_2\ldots a_k \quad\text{ and } \quad \b=b_1b_2\ldots b_k.
\end{equation*}
We call $(\a,\b)$ the {\sf endpoint notation} for $\pi$. For example, the endpoint notation for the partition in Figure \ref{fig:arcdiagram} is $(12357,25768)$. It is straightforward to check that partitions of $[n]$ are in bijection with pairs of vectors $(\a,\b)$ such that
\begin{itemize}
\item $\a$ and $\b$ have the same length (called the {\sf length} of the pair $(\a,\b)$),
\item $a_i<b_i$ for all $i$,
\item the entries of $\a$ are increasing, and
\item the entries of $\b$ are distinct.
\end{itemize}
In particular, the empty pair $(\emptyset,\emptyset)$ corresponds to the partition $\{\{1\},\{2\},\ldots,\{n\}\}$ and the longest pair $(12\ldots(n-1),23\ldots n)$ corresponds to the partition $\{[n]\}$. We will see that the endpoint notation is the best language for comparing Shi and Ish arrangements.

\subsection{Nonnesting partitions}
A partition $\pi$ of $[n]$ is called {\sf nonnesting} if it does not contain arcs $ij$ and $k\ell$ such that $i<k<\ell<j$ --- that is, no arc of $\pi$ ``nests'' inside another. The partition in Figure \ref{fig:arcdiagram} is {\em not} nonnesting (it is {\sf nesting}) because the arc $56$ nests inside the arc $37$. The number of nonnesting partitions of $[n]$ is famously given by the Catalan number $\frac{1}{n+1}\binom{2n}{n}$.

The property of nonnesting agrees well with the endpoint notation for partitions. That is, a partition $(\a,\b)$ is nonnesting if and only if its right endpoint vector $\b$ is {\em increasing}. In fact, the number of pairs of nesting arcs in $(\a,\b)$ is equal to the number of pairs $b_i>b_j$ such that $i<j$.

\subsection{$G$-partitions and $G$-Stirling numbers} Now we define a version of set partitions for any graph $G\subseteq\binom{[n]}{2}$: 
\begin{quote}
We say that a partition $\pi$ of $[n]$ is a {\sf $G$-partition} if all of its arcs are contained in the graph $G$.  The {\sf $G$-Stirling number} $\Stir(G,k)$ is the number of $G$-partitions with $k$ blocks.
\end{quote}
In particular, when $G$ is the complete graph $K_n=\binom{[n}{2}$ the $G$-partitions are unrestricted partitions of $[n]$ and the $G$-Stirling numbers are the classical Stirling numbers (of the second kind). 



\subsection{Connectivity}
Finally, we mention an auxiliary (nontrivial) result which we need for the proof of the Main Theorem. For $i\leq j$ we say that a partition $\pi$ of the set $\{i,i+1,\ldots,j\}$ is {\sf connected} if there does not exist $i\leq k<j$ such that $\pi$ refines the partition
\begin{equation*}
\{\{i,i+1,\dots,k\},\{k+1,\ldots,j-1,j\}\}.
\end{equation*}
Equivalently, $\pi$ is connected if its arc diagram has no holes when seen from space. The partition in Figure \ref{fig:arcdiagram} {\em is} connected. Moreover, a partition $\pi$ of $[n]$ has {\sf $d$ connected components} if there exist numbers $1<i_1,<\cdots <i_{d-1}<n$ such that $\pi$ refines the partition
\begin{equation*}
\{\{1,2,\ldots,i_1-1\},\{i_1,i_1+1,\ldots,i_2-1\},\ldots,\{i_{d-1},i_{d-1}+1,\ldots,n\}\}
\end{equation*}
and if its restriction to each block of this partition is connected. Equivalently, the arc diagram of a partition with $d$ connected components has $d-1$ holes when seen from space. For example, the partition $\{\{1,2\},\{3,5,7\},\{4,6\},\{8\}\}$ has $3$ connected components.

The second author has recently established an enumerative formula for nonnesting partitions (and other Catalan objects) that takes account of the type of the partition and its number of connected components. The prototype for this formula is the following theorem of Kreweras \cite[Theorem 4]{Kreweras}. Kreweras stated his formula in terms of noncrossing partitions; however, type-preserving bijections between noncrossing and nonnesting partitions have been observed by several authors.

\begin{lem}
\label{lem:kreweras}
Let $n>0$ and suppose that the sequence $(r_1,\ldots,r_n)$ of nonnegative integers satisfies $\sum_i i r_i = n$ and $\sum_i r_i = k$. The number of nonnesting partitions of $[n]$ with type $(r_1,\ldots,r_n)$ is
\begin{equation*}
\frac{ n! } {(n - k + 1)! r_1! r_2! \dots r_n!}.
\end{equation*}
\end{lem}

We will need the following formula of the second author \cite{RCat} in our proof of part (3) of the Main Theorem. The proof of this result is combinatorial and relies on the enumeration of words in certain monoids.


\begin{lem} \cite[Theorem 2.3, Part 2]{RCat}
\label{lem:rhoades}
Let $n > 0$ and suppose that the sequence $(r_1, \dots, r_n)$ of nonnegative integers satisfies $\sum_i i r_i = n$ and
$\sum_i r_i = k$.  Let $k \geq d$ and assume that $(r_1, \dots, r_n) \neq (n, 0, \dots, 0)$.  The number of nonnesting partitions of $[n]$ with type $(r_1, \dots, r_n)$ and $d$ connected components is
\begin{equation*}
\frac{ d (n - d - 1)! (k - 1)! } {(n - k - 1)! (k - d)! r_1! r_2! \dots r_n!}.
\end{equation*}
\end{lem}

When $k < d$, it is clear that there are no partitions of $[n]$ of type $(r_1, \dots, r_n)$ with $\sum_i r_i = k$ and $d$ connected components; there is a unique (nonnesting) partition of $[n]$ with type $(n, 0, \dots, 0)$ and it has $n$ connected components. 

We remark that the product formula in Lemma 2.3 was predicted from the formula \eqref{eq:candd} for Ish arrangements, {\bf not} by studying nonnesting partitions directly. This is one case in which the Ish arrangement acted as a ``toy model'' for other Catalan objects.

\section{Characteristic Polynomials}
In this section we explicitly compute the characteristic polynomials of $\Shi(G)$ and $\Ish(G)$ and observe that they are equal. The formula is expressed in terms of $G$-Stirling numbers $\Stir(G,k)$. Our tools are the finite field method of Crapo and Rota and the principle of inclusion-exclusion. Zaslavsky's theorem then implies that $\Shi(G)$ and $\Ish(G)$ have the same number of regions and the same number of relatively bounded regions (regions with one degree of freedom).

\subsection{The method}
Let $\A$ be a finite hyperplane arrangement in $\RR^n$ and suppose that the defining equations for hyperplanes in $\A$ have coefficients in $\ZZ$. Then the finite field method of Crapo and Rota \cite{CrapoRota} is a useful way to compute the characteristic polynomial of $\A$ without having to know its intersection poset. Let $p\in\ZZ$ be prime and consider a hyperplane $H\subseteq\RR^n$ with fixed defining equation $a_1x_1+\cdots +a_nx_n=b$, where $a_i,b\in\ZZ$. Then we define the following subset $H_p$ of the {\em finite} vector space $\FF_p^n$ by reducing the coefficients of $H$ modulo $p$:
\begin{equation*}
H_p :=\{ (x_1,x_2,\ldots,x_n)\in \FF_p^n : a_1x_1+\cdots + a_nx_n = b\}.
\end{equation*}
Observe that $H_p$ may not be a hyperplane in $\FF_p^n$ when $p$ is small, and that $H_p$ in general depends on the defining equation chosen. However: If $p$ is large enough then each $H_p$ {\em is} a hyperplane in $\FF_p^n$ and the characteristic polynomial of $\A$ has a nice relationship to the {\sf reduced hyperplane arrangement} $\A_p := \{ H_p : H\in\A\}$ in $\FF_p^n$.

\begin{thm} \cite{CrapoRota}  Let $p\in\ZZ$ be a large prime, and let $\A$ be a finite collection of hyperplanes in $\RR^n$ whose hyperplanes have defining equations with coefficients in $\ZZ$. Then the characteristic polynomial of $\A$ satisfies
\begin{equation*}
\chi_\A(p) = \#\left( \FF_p^n - \cup_{H\in\A} H_p\right).
\end{equation*}
That is, $\chi_\A(p)$ counts the number of points in the complement of the reduced arrangement $\A_p$ in the finite vector space $\FF_p^n$.
\end{thm}

\subsection{The calculation}
Now we use the finite field method to compute the characteristic polynomials of $\Shi(G)$ and $\Ish(G)$. We observe that they are equal.

\begin{thm}
\label{th:charpoly}
Let $G\subseteq\binom{[n]}{2}$ be a graph on $n$ vertices. The characteristic polynomials of the deleted Shi and Ish arrangement are given by:
\begin{equation*}
\chi_{\Shi(G)}(p)=\chi_{\Ish(G)}(p)=p \sum_{k = 0}^{n-1} (-1)^k \Stir(G,n-k) \frac{(p-k-1)!}{(p-n)!}.
\end{equation*}
\end{thm}

\begin{proof}
Let $p\in\ZZ$ be a large prime. We will show that the reduced complements $\FF_p^n-\Shi(G)_p$ and $\FF_p^n-\Ish(G)_p$ (forgive the abuse of notation) contain the same number of points, counted by the above formula.

To do this, we identify $\{0,1,\ldots,p-1\}=\FF_p$ with the vertices of a regular $p$-gon, ordered clockwise. (That is, $i+1$ is just clockwise of $i$.) Then a vector $v=(v_1,\ldots,v_n)\in\FF_p^n$ is a labeling of the vertices: if $v_i=j$ then we place the label $v_i$ on the vertex $j$. Note that $v\in\FF_p^n$ is in the complement of the (reduced) Coxeter arrangement $\Cox(n)_p$ precisely when $v_i-v_j\neq 0$ for all $1\leq i<j\leq n$. That is, the points of $\FF_p^n-\Cox(n)_p$ correspond to {\em injective} labelings $\{v_1,\ldots,v_n\}\hookrightarrow \FF_p$. The complements of $\Shi(G)_p$ and $\Ish(G)_p$ are both contained in $\FF_p^n-\Cox(n)_p$, so we must count certain kinds of injective labelings.

First we deal with $\Shi(G)_p$. For any set of edges $S\subseteq G$ let $f(S)$ denote the number of vectors $v\in\FF_p^n-\Shi(G)_p$ such that $v_i-v_j=1$ for all edges $ij\in S$ (this notation implies $i<j$). By the principle of inclusion-exclusion (see for example \cite[Chapter 2]{StanEC1}) we observe that the number of points in $\FF_p^n-\Shi(G)_p$ is equal to
\begin{equation}
\label{eq:finclusion}
\sum_{S\subseteq G} (-1)^{|S|}f(S).
\end{equation}
Now suppose that $S$ contains edges $ij$ and $i\ell$ with the same left endpoint. The conditions $v_i-v_j=1$ and $v_i-v_\ell=1$ imply that $v_j=v_\ell$ which cannot be satisfied on $\FF_p^n-\Cox(n)_p$, hence $f(S)=0$. Similarly $f(S)=0$ whenever $S$ contains two edges with the same right endpoint. That is, the sets $S$ that contribute to the sum \eqref{eq:finclusion} are precisely the arc sets of $G$-partitions.

Let $S\subseteq G$ correspond to a $G$-partition with $n-k$ blocks (that is, $|S|=k$). To compute $f(S)$ note that the conditions $v_i-v_j=1$ for all $ij\in S$ imply that the $p$-gon $\FF_p$ gets labeled by $n-k-1$ (given) contiguous strings of labels with spaces between. There are $(n-k-1)!$ ways to cyclically permute the strings; there are $\binom{p-k-1}{n-k-1}$ ways to place $p-n$ empty spaces between the strings; and there are $p$ ways to choose the ``origin'' (the location of $0$). Hence:
\begin{equation}
\label{eq:fformula}
f(S)=\frac{p\,(p-k-1)!}{(p-n)!}.
\end{equation}
Combining \eqref{eq:finclusion} and \eqref{eq:fformula} with the finite field method gives the desired formula.

We use a parallel argument to deal with $\Ish(G)_p$. For any set of edges $S\subseteq G$ let $g(S)$ be the number of vectors $v\in\FF_p^n-\Ish(G)_p$ such that $v_1-v_j=i$ for all $ij\in S$. As above, the points of $\FF_p-\Ish(G)_p$ are counted by
\begin{equation}
\sum_{S\subseteq G} (-1)^{|S|}g(S),
\end{equation}
and one can check that $g(S)=0$ unless $S$ is the arc set of a $G$-partition. We let $S$ correspond to a $G$-partition with $n-k$ blocks (that is, $|S|=k$) and compute $g(S)$ as follows. First choose $v_1$ in $p$ ways. Then for each $ij\in S$ the condition $v_1-v_j=i$ uniquely determines the value of $v_j$. The remaining $n-k-1$ labels must be placed injectively in the remaining $p-k-1$ positions and there are $(p-k-1)(p-k-2)\cdots (p-n+1)$ ways to do this. Thus we get the desired formula:
\begin{equation}
g(S)=\frac{p\,(p-k-1)!}{(p-n)!}.
\end{equation} 
\end{proof}

Notice that the counting argument for computing $\chi_{\Ish(G)}(p)$ was more straightforward than the argument for $\chi_{\Shi(G)}(p)$. This again agrees with our observation that Ish is a toy model for Shi. It is somewhat surprising that the two inclusion-exclusion arguments result in the same expression. It may be interesting to find a direct bijection between the points of the complements $\FF_p^n-\Shi(G)_p$ and $\FF_p^n-\Ish(G)_p$.

\subsection{Remarks} A simplified version of the above argument shows that the characteristic polynomials of $\Shi(n)$ and $\Ish(n)$ (the case of the complete graph) are both equal to $p\,(p-n)^{n-1}$. This result was obtained earlier by Headley \cite{Headley} and Athanasiadis \cite{Athanasiadis} (for the Shi arrangement) and by the first author \cite[Theorem 1]{Armqt} (for the Ish arrangement). Moreover, Athanasiadis described a special family of graphs $G$ for which the characteristic polynomial of $\Shi(G)$ splits. His result \cite[Theorem 2.2]{AFree} together with Theorem \ref{th:charpoly} implies the following.


\begin{cor}
Suppose the graph $G\subseteq\binom{[n]}{2}$ has the following property: if $i<j<k$ and $ij\in G$, then $ik\in G$. Then we have
\begin{equation*}
\chi_{\Shi(G)}(p)=\chi_{\Ish(G)}(p)=p\prod_{i=1}^{n-1} (p-d_i-i),
\end{equation*}
where $d_i:=\#\{j : ij\in G\}$ is the outdegree of vertex $i$ in $G$.
\end{cor}
In the same paper, Athanasiadis showed that the arrangements $\Shi(G)$ of the Corollary are {\sf free} in the sense of Terao \cite{Terao} (see also
\cite{OrlikTerao}). 
This is an open problem for the corresponding Ish arrangements $\Ish(G)$. 


We remark that the characteristic polynomial of an arrangement allows us to count certain kinds of regions. Some notation: Let $\A$ be a finite collection of hyperplanes in $\RR^n$ and suppose that the normals to the hyperplanes span a space $V\subseteq\RR^n$ of dimension $r$. This $r$ is called the {\sf rank} of the arrangement. If $r<n$ then the arrangement $\A$ has no bounded regions; in this case we say that a region of $\A$ is {\sf relatively bounded} if its intersection with $V$ is bounded. The following is a classic theorem of Zaslavsky.

\begin{thm} \cite{Zaslavsky} Let $\A$ be a hyperplane arrangement in $\RR^n$ with rank $r$. Then:
\begin{itemize}
\item The number of regions of $\A$ is $(-1)^n \chi_\A(-1)$;
\item The number of relatively bounded regions of $\A$ is $(-1)^r \chi_\A(1)$.
\end{itemize}
\end{thm}

\begin{cor}
The arrangements $\Shi(G)$ and $\Ish(G)$ have the same number of regions and the same number of relatively bounded regions.
\end{cor}

Observe that the normals to either $\Shi(G)$ or $\Ish(G)$ span the hyperplane $x_1+x_2+\cdots + x_n=0$. Hence each of these arrangements has rank $n-1$. It follows that neither arrangement has bounded regions and its relatively bounded regions have one degree of freedom. In the case of the complete graph, we find that the arrangements $\Shi(n)$ and $\Ish(n)$ both have $(n+1)^{n-1}$ regions and $(n-1)^{n-1}$ regions with one degree of freedom.

The fact that the Shi arrangement $\Shi(n)$ has $(n+1)^{n-1}$ regions was first proved by Jian-Yi Shi (see \cite{Shi}). This beautiful result has motivated more than a few research papers since 1985 (including the present one).

\section{Labeling the regions} 
Fix a graph $G\subseteq\binom{[n]}{2}$. In this section we devise combinatorial labels for the regions of the deleted arrangements $\Shi(G)$ and $\Ish(G)$; we call these labels {\sf Shi ceiling diagrams} and {\sf Ish ceiling diagrams}, respectively. (Something like ``Shi floor diagrams'' appeared earlier in Athanasiadis and Linusson \cite{ALShi}.) Essentially, each diagram encodes the ceilings of a given region, from which we can easily determine its recession cone.

\subsection{Shi ceiling diagrams} Recall that the regions (cones) of the Coxeter arrangement $\Cox(n)$ correspond to elements of the symmetric group $\frak{S}(n)$. If $C\subseteq\RR^n$ is the {\sf dominant cone} --- defined by the coordinate inequalities $x_1>x_2>\cdots > x_n$ --- then the collection of regions of $\Cox(n)$ is $\{wC : w\in\frak{S}(n)$\}, where $wC$ is defined by the coordinate inequalities
\begin{equation}
\label{eq:shiwcone}
x_{w(1)}>x_{w(2)}>\cdots >x_{w(n)}.
\end{equation}

Now let $R$ be a region of the deleted Shi arrangement $\Shi(G)$. Since $\Cox(n)\subseteq\Shi(G)$, $R$ is contained in some cone $wC$. In this case, what are the possible ceilings of $R$? We note that the hyperplanes of $\Shi(G)$ that intersect $wC$ are precisely
\begin{equation*}
\Phi^+(G,w):=\{ x_{w(i)}-x_{w(j)} = 1 : ij\in G \text{ and } w(i)<w(j) \}.
\end{equation*}
(We can think of these as the {\em non-inversions} of $w$ contained in $G$.) Furthermore, suppose that the region $R$ is ``below'' some hyperplane $x_{w(i)}-x_{w(j)}=1$ --- that is, suppose that each $v\in R$ satisfies $v_{w(i)}-v_{w(j)}<i$. Then, considering \eqref{eq:shiwcone}, $R$ is also below any hyperplane of the form $x_{w(i')}-x_{w(j')}=1$ such that
\begin{equation}
\label{eq:rootposet}
w(i)\leq w(i')<w(j')\leq w(j').
\end{equation}
That is, if we declare a partial order on $\Phi^+(G,w)$ by saying that $x_{w(i')}-x_{w(j')}=1$ is ``less than'' $x_{w(i)}-x_{w(j)}=1$ when condition \eqref{eq:rootposet} holds, then the collection of $\Shi(G)$-hyperplanes above $R$ forms a down-closed set.

\begin{thm}
There is a bijection between regions of $\Shi(G)$ in the cone $wC$ and {\sf order ideals} (down-closed sets) in the poset $\Phi^+(G,w)$. This map sends a region $R$ to the set of hyperplanes in $\Shi(G)$ that are ``above'' $R$ (contain $R$ and the origin in the same half space). The {\em maximal} elements of the ideal are the ceilings of $R$.
\end{thm}

\begin{proof}
Let $R$ be a region of $\Shi(G)$ contained in $wC$. We showed above that the collection of hyperplanes above $R$ is an order ideal in $\Phi^+(G,w)$. The map is injective since these hyperplanes uniquely determine $R$. Observe that the ceilings of $R$ are the elements of the ideal that may be individually removed to obtain another ideal, and these are precisely the maximal elements. We refer to Athanasiadis and Linusson \cite{ALShi} for the proof that every ideal corresponds to a non-empty region.
\end{proof}

To express this combinatorially, {\bf we note that order ideals in $\Phi^+(G,w)$ are equivalent to nonnesting $G$-partitions whose blocks are ``increasing'' with respect to $w$.} Indeed, there is a bijection between ideals and {\sf antichains} (sets of pairwise-incomparable elements), since an ideal is uniquely determined by its antichain of maximal elements. By sending the hyperplane $x_{w(i)}-x_{w(j)}=1$ to the arc $ij$, each antichain in $\Phi^+(G,w)$ corresponds to a $G$-partition of $[n]$ whose arcs $i<j$ satisfy $w(i)<w(j)$. Finally, note that two arcs nest if and only if they are comparable in the poset $\Phi^+(G,w)$.

Following these remarks, we draw a diagram for each region of $\Shi(G)$.
\begin{defn}
Let $R$ be a region of $\Shi(G)$ contained in the cone $wC$. We associate $R$ with the pair $(w,\pi)$ where $\pi$ is an order ideal in the poset $\Phi^+(G,w)$ of non-inversions of $w$ contained in $G$. Equivalently, $\pi$ is a nonnesting $G$-partition whose blocks are increasing with respect to $w$. We draw $(w,\pi)$ by placing the arc diagram for $\pi$ above the numbers $w(1),\ldots,w(n)$, and we call this the {\sf Shi ceiling diagram} of $R$.
\end{defn}

\begin{figure}
\begin{center}
\includegraphics[scale=1]{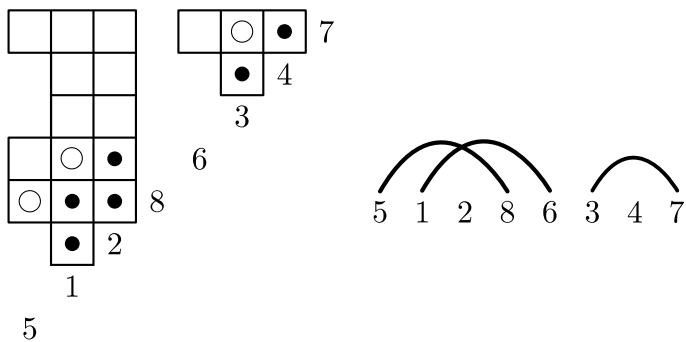}
\end{center}
\caption{An order ideal (left) and a  Shi ceiling diagram (right)}
\label{fig:shiceilingdiagram}
\end{figure}

For example, let $K_8=\binom{[8]}{2}$ be the complete graph on $8$ vertices and consider the permutation $w=51286347\in\frak{S}(8)$. Figure \ref{fig:shiceilingdiagram} displays the ideal in $\Phi^+(K_8,w)$ (left) and the ceiling diagram (right) corresponding to a region $R$ of $\Shi(K_8)=\Shi(8)$ contained in the cone $wC$. The squares are elements of the poset $\Phi^+(K_8,w)$ and the circles are elements of the ideal (closed to the right and down). The hollow circles (maximal elements) indicate the ceilings of the region: $x_5-x_8=1$, $x_1-x_6=1$, $x_3-x_7=1$. The corresponding nonnesting partition is $\pi=\{\{1,4\},\{2,5\},\{3\},\{6,8\},\{7\}\}$.

\begin{figure}
\begin{center}
\includegraphics[scale=1]{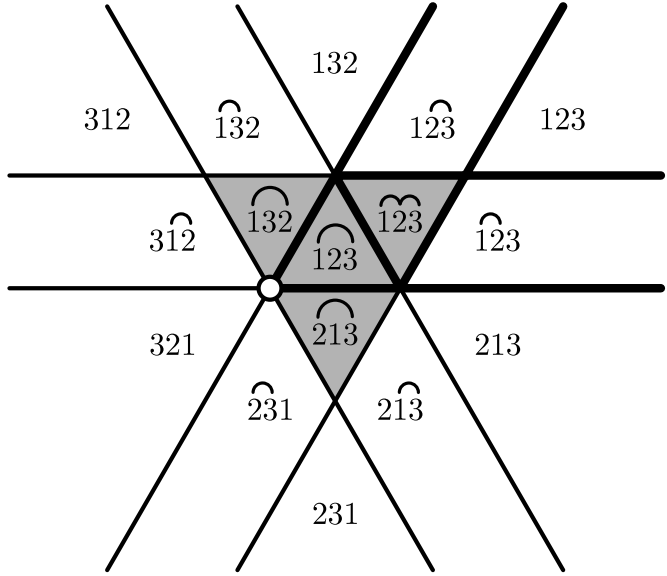}
\end{center}
\caption{The Shi arrangement $\Shi(3)$ labeled by ceiling diagrams}
\label{fig:shiceilings}
\end{figure}

Figure \ref{fig:shiceilings} displays the whole arrangement $\Shi(3)$ with its regions labeled by ceiling diagrams. Observe that we can read the degrees of freedom from the ceiling diagram $(w,\pi)$: the corresponding region has $d$ degrees of freedom if and only if the nonnesting partition $\pi$ has $d$ connected components. This is a general phenomenon.

\begin{lem}
Let $R$ be a region of $\Shi(G)$ with ceiling diagram $(w,\pi)$. This region has $d$ degrees of freedom if and only if the nonnesting partition $\pi$ of $[n]$ has $d$ connected components.
\end{lem}

\begin{proof}
Suppose that $\pi$ has $d$ connected components. That is, there exist $1<i_1<\cdots <i_{d-1}<n$ such that $\pi$ refines the partition
\begin{equation*}
\{ \{1, 2, \dots, i_1 - 1 \}, \{i_1, i_1 + 1, \dots, i_2 - 1\} ,\dots, \{i_{d-1}, i_{d-1} + 1, \dots, n \} \},
\end{equation*}
and its restriction to any block of this partition is connected. We compute the recession cone of $\Rec(R)\subseteq\RR^n$ of $R$ as follows.

Consider $v=(v_1,\ldots,v_n)\in\Rec(R)$. Since $R$ is in the cone $wC$ we must have $v_{w(1)}\geq v_{w(2)}\geq \cdots \geq v_{w(n)}$. Moreover, if $ij\in G$ with $i$ and $j$ in the same block of $\pi$ then the coordinate inequality $x_{w(i)}-x_{w(j)}<1$ holds on $R$ and we must have $v_{w(i)}=v_{w(i+1)}=\cdots=v_{w(j)}$. Since these are the only constraints on $v$, we conclude that the recession cone $\Rec(R)$ consists of all vectors of the form $w\cdot (a_1,a_2,\ldots,a_n)$, where $a_1\geq \cdots \geq a_n$ and where $a_i=a_j$ if $i$ and $j$ are in the same connected component of $\pi$. The dimension of the cone is therefore $d$.
\end{proof}

For example, consider the ceiling diagram $(w,\pi)$ in Figure \ref{fig:shiceilingdiagram} and the corresponding region $R$ of $\Shi(8)$. The connected components of $\pi$ are $\{1,2,3,4,5\},\{6,7,8\}$ and their images under $w$ are $\{5,1,2,8,5\},\{3,4,7\}$. Hence the recession cone $\Rec(R)$ consists of all vectors of the form $(a,a,b,b,a,a,b,a)\in\RR^8$ with $a\geq b$, and it has dimension $2$.

\subsection{Ish ceiling diagrams}
In order to compare the two arrangements, we now define an Ish analogue of Shi ceiling diagrams.

Since $\Cox(n)\subseteq \Ish(G)$, each region $R$ of $\Ish(G)$ is contained in $wC$ for some permutation $w\in\frak{S}(n)$, in which case each vector $v=(v_1,\ldots,v_n)\in R$ satisfies
\begin{equation}
\label{eq:wcone}
v_{w(1)}> v_{w(2)}> \cdots > v_{w(n)}.
\end{equation}
Which $\Ish(G)$-hyperplanes are the possible ceilings of this region? If the hyperplane $x_1-x_j=i$ intersects the cone $wC$ it must be true that $x_1>x_j$ on $wC$ (since $i$ is positive). Considering \eqref{eq:wcone}, this means that $j$ must occur to the right of $1$ in the list $w(1),\ldots,w(n)$ --- that is, we must have $w^{-1}(1)<w^{-1}(j)$. We conclude that the Ish hyperplanes that intersect the cone $wC$ are precisely
\begin{equation*}
\Psi^+(G,w):=\{ x_1-x_j=i : ij\in G \text{ and } w^{-1}(1)<w^{-1}(j)\}.
\end{equation*}
Now let $R$ be a region of $\Ish(G)$ in the cone $wC$ and suppose that $R$ is below $x_1-x_j=i$ --- that is, each $v\in R$ satisfies $v_1-v_j<i$. Then it is easy to check that $R$ is also below the hyperplane $x_1-x_{j'}=i'$, where
\begin{equation}
\label{eq:ishroots}
\text{ either } i<i' \text{ or } w^{-1}(j')<w^{-1}(j).
\end{equation}
By analogy with the Shi case, we define a partial order on $\Psi^+(G,w)$ by declaring that the hyperplane $x_1-x_j=i$ is ``less than'' the hyperplane $x_1-x_{j'}=i'$ whenever \eqref{eq:ishroots} holds. This leads to a useful characterization of $\Ish(G)$ regions.

\begin{thm}
There is a bijection between regions of $\Ish(G)$ in the cone $wC$ and {\sf order filters} (up-closed sets) in the poset $\Psi^+(G,w)$. This map sends a region $R$ to the set of hyperplanes in $\Ish(G)$ that are ``above'' $R$ (contain $R$ and the origin in the same half space). The {\em minimal} elements of the filter are the ceilings of $R$.
\end{thm}

\begin{proof}
Let $R$ be a region of $\Ish(G)$ in the cone $wC$. By the above remarks we know that the collection of $\Ish(G)$-hyperplanes above $R$ is an order filter in $\Psi^+(G,w)$. These hyperplanes 
together with the fact that $R$ lies in $wC$
uniquely determine $R$, so the map is injective. The ceilings of $R$ are precisely the elements of this filter that may be individually removed to obtain another filter --- that is, they are the minimal elements.

To show that the map is surjective, we must show that each filter in $\Psi^+(G,w)$ corresponds to a non-empty region of $\Ish(G)$.  Let 
$F \subseteq \Psi^+(G,w)$ be an order filter and let $A \subseteq F$ be its set of 
minimal elements.  For $1 \leq i \leq n$ define 
\begin{equation*}
z_i = -\mathrm{max}\{ j \,:\, \text{$x_1 - x_k = j \in A$ and $w^{-1}(k) \leq w^{-1}(i)$} \},
\end{equation*}
where we adopt the convention that $\mathrm{max}(\emptyset) = 0$.  One may check that the point $(z_{w(1)}, \dots, z_{w(n)}) \in \RR^n$ lies on the boundary of a region of $\Ish(G)$ which maps to the filter $F$. (Alternatively, note that Theorems \ref{th:charpoly} and \ref{thm:candd} imply, respectively, that the number of regions of $\Ish(G)$ and the number of filters in $\Psi^+(G,w)$ (summed over $w$) are both equal to
\begin{equation*}
\sum_{k=0}^{n-1} \Stir(G,n-k) \frac{n!}{(k+1)!}.
\end{equation*}
Hence any injective map between then must be surjective.)
\end{proof}

\begin{figure}
\begin{center}
\includegraphics[scale=1]{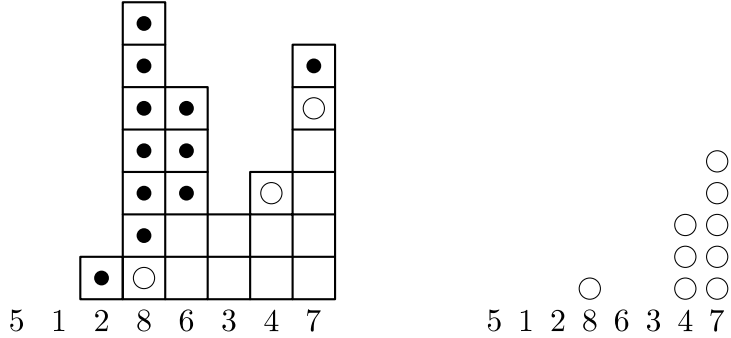}
\end{center}
\caption{An order filter (left) and an Ish ceiling diagram (right)}
\label{fig:ishceilingdiagram}
\end{figure}

It is convenient to express this situation with a picture. Given $w\in\frak{S}(n)$ we draw $w(1),w(2),\ldots,w(n)$ on a line. For each $j$ to the right of $1$ we draw $j-1$ boxes above the symbol $j$. If we identify the $i$th box above $j$ with the hyperplane $x_1-x_j=i$ then the collection of boxes is exactly $\Psi^+(G,w)$ (we erase the boxes that are not in $G$); the partial order on boxes increases up and to the left.

\begin{figure}
\begin{center}
\includegraphics[scale=1]{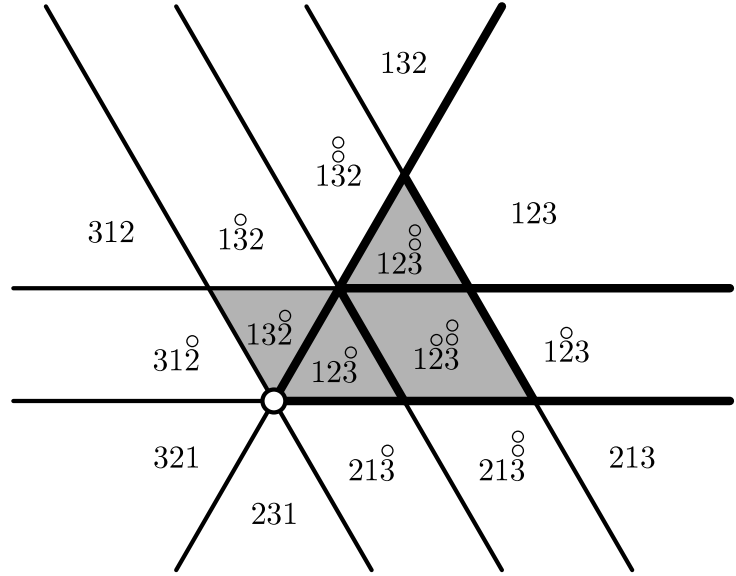}
\end{center}
\caption{The Ish arrangement $\Ish(3)$ labeled by ceiling diagrams}
\label{fig:ishceilings}
\end{figure}

In Figure \ref{fig:ishceilingdiagram} (left side) we have drawn the poset $\Psi^+(G,w)$ for the complete graph $G=K_8$ and the permutation $w=51286347\in\frak{S}(8)$. The circles (closed up and to the left) indicate an order filter in this poset. This filter defines a region $R$ of $\Ish(K_8)=\Ish(8)$ in the cone $wC$ and its ceilings are the antichain of minimal elements (hollow circles): $x_1-x_8=1$, $x_1-x_4=3$, $x_1-x_7=5$.  To simplify the diagram further (right side), we just draw $i$ hollow circles above the symbol $j$ for each ceiling $x_1-x_j=i$. This is the {\sf Ish ceiling diagram} of the region. We will encode it with the pair $(w,\varepsilon)$ where $\varepsilon_i$ is the number of circles above the symbol $w(i)$. For the example in Figure \ref{fig:ishceilingdiagram} we have
\begin{equation*}
(w,\varepsilon)=(51286347, (0,0,0,1,0,0,3,5)).
\end{equation*}
Figure \ref{fig:ishceilings} displays the full arrangement $\Ish(3)$ with regions labeled by ceiling diagrams.

In order to count regions later, here is a purely combinatorial characterization of Ish ceiling diagrams.

\begin{defn}
\label{defn:ishceilingdiagram}
Let $G\subseteq\binom{[n]}{2}$ be a graph and consider a permutation $w\in\frak{S}(n)$. We call the pair $(w,\varepsilon)$ an {\sf Ish ceiling diagram} if the vector $\varepsilon=(\varepsilon_1,\ldots,\varepsilon_n)$ satisfies:
\begin{itemize}
\item $0\leq \varepsilon_i<w(i)$;
\item $\varepsilon_i=0$ unless $w^{-1}(1)<w^{-1}(i)$;
\item If $\varepsilon_i>0$ then $\varepsilon_i<w(i)$ is an edge in $G$;
\item the nonzero entries of $\varepsilon$ strictly increase.
\end{itemize}
We will draw the pair $(w,\varepsilon)$ by placing $w(1),\ldots,w(n)$ on a line and drawing $\varepsilon_i$ circles above $w(i)$. By the above remarks, the pair $(w,\varepsilon)$ corresponds to a unique region of $\Ish(G)$ with a ceiling $x_1-x_{w(i)}=i$ for each $\varepsilon_i\neq 0$.
\end{defn}

Finally, we can read the recession cone of a region directly from its Ish ceiling diagram.

\begin{lem}
Let $R$ be a region of $\Ish(G)$ in the cone $wC$ with ceiling diagram $(w,\varepsilon)$. If $k$ is the {\em maximum} index such that $\varepsilon_k\neq 0$ (or $k=w^{-1}(1)$ if $\varepsilon$ is the zero vector) then $R$ has $n-k+w^{-1}(1)$ degrees of freedom. In particular, the region $R$ is relatively bounded (has $1$ degree of freedom) if and only if $w(1)=1$ and $\varepsilon_n\neq 0$.
\end{lem}

\begin{proof}
Consider $v=(v_1,\ldots,v_n)\in\Rec(R)$. Since $R$ is in the cone $wC$ we must have $v_{w(1)}\geq v_{w(2)}\geq \cdots \geq v_{w(n)}$. If $\varepsilon_j\neq 0$ then we must also have $v_1-v_{w(j)}<\varepsilon_j$, which implies that $v_{w^{-1}(1)}=v_{w^{-1}(1)+1}=\cdots =v_{w^{-1}(j)}$. Since these are the only constraints on $v$, we conclude that the recession cone $\Rec(R)$ consists of all vectors of the form $w\cdot (a_1,a_2,\ldots,a_n)$ with $a_1\geq \cdots \geq a_n$ and $a_{w^{-1}(1)}=a_{w^{-1}(1)+1}=\cdots =a_k$. The dimension of the cone is therefore $n-(k-w^{-1}(1))=n-k+w^{-1}(1)$.
\end{proof}

For example, consider the Ish ceiling diagram $(w,\varepsilon)$ in Figure \ref{fig:ishceilingdiagram} and the corresponding region $R$ of $\Ish(8)$. In this case we have $w(1)=5$ and $k=8$ is the largest index such that $\varepsilon_k\neq 0$. Hence the recession cone $\Rec(R)$ consists of all vectors $(a,a,a,a,b,a,a,a)\in\RR^8$ with $a\geq b$, and it has dimension $2$.

\subsection{A bijection between dominant regions}
The Shi and Ish ceiling diagrams immediately give us a bijection between {\bf dominant} regions of $\Shi(G)$ and $\Ish(G)$ with the same number of ceilings. This bijection does {\bf not} preserve degrees of freedom because it can't: in general $\Shi(G)$ and $\Ish(G)$ have different numbers of relatively bounded dominant regions. For example, $\Shi(3)$ has $2$ (see Figure \ref{fig:shiceilings}) and $\Ish(3)$ has $3$ (see Figure \ref{fig:ishceilings}).

\begin{thm}
\label{thm:dominant}
Consider a graph $G\subseteq\binom{[n]}{2}$ and an integer $c$. The deleted arrangements $\Shi(G)$ and $\Ish(G)$ have the same number of {\bf dominant} regions with $c$ ceilings.
\end{thm}

\begin{proof}
This is essentially a picture proof. For the identity permutation $w={\bf 1}$ we observe that the posets $\Phi^+(G,{\bf 1})$ and $\Psi^+(G,{\bf 1})$ look exactly the same, except that one is reflected in a line of slope $1$. For example, here are the posets corresponding to the graph $G=\binom{[8]}{2}-\{14,34,48,58\}$; Shi on the left, Ish on the right:
\begin{center}
\includegraphics[scale=1]{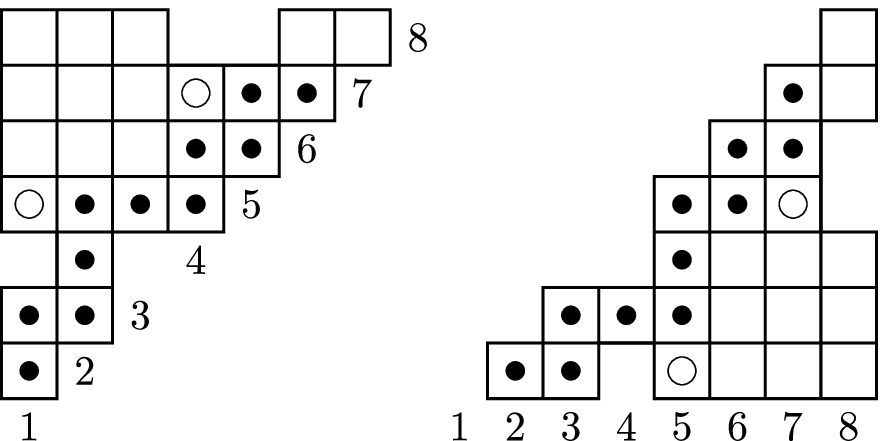}
\end{center}
This reflection is an order-reversing bijection between $\Phi^+(G,{\bf 1})$ and $\Psi^+(G,{\bf 1})$. Hence it induces a bijection between {\bf ideals in $\Phi^+(G,{\bf 1})$ with $c$ maximal elements} (dominant $\Shi(G)$-regions with $c$ ceilings) and {\bf filters in $\Psi^+(G,{\bf 1})$ with $c$ minimal elements} (dominant $\Ish(G)$-regions with $c$ ceilings).
\end{proof}

The number of dominant regions with $c$ ceilings equals the Narayana number $\frac{1}{n}\binom{n}{c}\binom{n-1}{c}$ when $G$ is the complete graph, and equals the binomial coefficient $\binom{n-1}{c}$ when $G$ is the chain $\{12,23,\ldots,(n-1)n\}$. We do not know a closed formula for general $G$.

Note that the bijection in Theorem \ref{thm:dominant} does {\bf not} extend to other cones $wC$, since in general the posets $\Phi^+(G,w)$ and $\Psi^+(G,w)$ look very different. Indeed, $\Shi(G)$ and $\Ish(G)$ do {\bf not} have the same number of regions in a given cone $wC$. (Consider Figures \ref{fig:shiceilings} and \ref{fig:ishceilings} with the permutation $w=132$.)

However, we gain something by summing over the cones $wC$. Not only do $\Shi(G)$ and $\Ish(G)$ have the same number of (unrestricted) regions with $c$ ceilings, they have the same number of regions with $c$ ceilings and $d$ degrees of freedom. We prove this in the next section using a {\em non-bijective} method.

\section{Counting the regions}
In this section we introduce a partition-valued statistic on the regions of $\Shi(G)$ and $\Ish(G)$, and in each case we call this the {\sf ceiling partition} of the region. (This concept is new even for the Shi arrangement.) It turns out that $\Shi(G)$ and $\Ish(G)$ have the same number of regions $R$ with a given ceiling partition $\pi$; moreover, when the partition $\pi$ has $k$ blocks (i.e. $R$ has $n-k$ ceilings), this number has a beautiful formula: $n!/(n-k+1)!$. The partition $\pi$ does {\bf not} determine the degrees of freedom of $R$. However, we still have a nice formula: The number of regions of $\Shi(G)$ or $\Ish(G)$ with ceiling partition $\pi$ (with $k$ blocks) and $d$ degrees of freedom equals
\begin{equation*}
\frac{d(n-d-1)!(k-1)!}{(n-k-1)!(k-d)!}.
\end{equation*}
This completes the proof of the Main Theorem. At the end we make comments and suggestions for future research.

\subsection{Ceiling partitions}
To each region $R$ of $\Shi(G)$ or $\Ish(G)$ we associate a partition of the set $[n]$, called its {\sf ceiling partition}. We note that this partition may be {\em nesting}, and in general every partition of $[n]$ will occur. The ceiling partition is determined by the ceilings of $R$ and the cone $wC$ in which $R$ occurs; thus we can read it from the ceiling diagram. We will see that the correct language for ceiling partitions is the {\sf endpoint notation} $(\a,\b)$, discussed in Section \ref{sec:partitions}.

\begin{defn}\hspace{.1in}
\begin{enumerate}
\item Let $R$ be a region of $\Shi(G)$ with ceiling diagram $(w,\pi)$. Then the {\sf ceiling partition} of $R$ is $w\cdot \pi$ ($w$ acting on $\pi$). That is, the ceiling partition has $w(i)$ and $w(j)$ in a block whenever $i$ and $j$ are in a block of $\pi$. For example, the region in Figure \ref{fig:shiceilingdiagram} has ceiling partition $\{\{1,6\},\{2\},\{3,7\},\{4\},\{5,8\}\}$, with endpoint notation $(\a,\b)=(135,678)$. Note that the ceiling partition $(\a,\b)$ has $c$ arcs if and only if $R$ has $c$ ceilings.

\medskip

\item Let $R$ be a region of $\Ish(G)$ with ceiling diagram $(w,\varepsilon)$. We define a pair of vectors $(\a,\b)$ such that $a_i$ is the $i$th nonzero entry of $\varepsilon$, which occurs in position $w^{-1}(b_i)$. The conditions of Definition \ref{defn:ishceilingdiagram} guarantee that $(\a,\b)$ is the endpoint notation for a partition, which we call the {\sf ceiling partition} of $R$. For example, the region shown in Figure \ref{fig:ishceilingdiagram} has ceiling partition $(\a,\b)=(135,847)$ since there is one circle above $8$, three above $4$, and five above $7$. Again, the ceiling partition has $c$ arcs if and only if $R$ has $c$ ceilings.
\end{enumerate}
\end{defn}

\subsection{Counting Shi and Ish regions}
Let $c$ and $d$ be integers. We separately count the regions of $\Shi(G)$ and $\Ish(G)$ with $c$ ceilings and $d$ degrees of freedom, and observe that they are the same. This completes the proof of the Main Theorem.

\begin{thm}
\label{thm:candd}
Fix a graph $G\subseteq\binom{[n]}{2}$ and let $\A$ be either $\Shi(G)$ or $\Ish(G)$. Let  $(\a,\b)$ be a partition of $[n]$ with $k$ blocks ($n-k$ arcs) and consider an integer $1\leq d\leq k$. There exists a region of $\A$ with ceiling partition $(\a,\b)$ if and only if we have $a_ib_i\in G$ for all $i$, in which case:
\begin{enumerate}
\item The number of regions of $\A$ with ceiling partition $(\a,\b)$ is
\begin{equation*}
\frac{n!}{(n-k+1)!}.
\end{equation*}
\item The number of regions of $\A$ with ceiling partition $(\a,\b)$ and $d$ degrees of freedom is
\begin{equation}
\label{eq:candd}
\frac{d(n-d-1)!(k-1)!}{(n-k-1)!(k-d)!}.
\end{equation}
\end{enumerate}
To obtain the number of regions with $c$ ceilings and $d$ degrees of freedom, sum \eqref{eq:candd} over $G$-partitions $(\a,\b)$ with $k=n-c$ blocks.
\end{thm}
\begin{proof}
First we deal with $\A=\Shi(G)$. Recall that a Shi ceiling diagram $(w,\pi)$ is a nonnesting partition $\pi$ whose blocks are increasing with respect to the permutation $w$. Thus, to create a ceiling diagram (region) with ceiling partition $(\a,\b)$, we must first choose a nonnesting partition $\pi_0$ with the same block sizes as $(\a,\b)$ and then put the labels from each block of $(\a,\b)$ (increasingly) in a block of $\pi_0$. So suppose that $(\a,\b)$ has $r_i$ blocks of size $i$. By Lemma \ref{lem:kreweras} there are
\begin{equation*}
\frac{n!}{(n-k+1)!r_1!r_2!\cdots r_n!}
\end{equation*}
ways
to choose $\pi_0$. Then, there are $r_1!r_2!\cdots r_n!$ ways to map each block of $(\a,\b)$ to a block of $\pi_0$ with the same size. This proves (1). To prove (2), note that the region $(w,\pi)$ has $d$ degrees of freedom if and only if the nonnesting partition $\pi$ has $d$ connected components. Apply the same argument as above, but use Lemma \ref{lem:rhoades}.

Next we deal with $\A=\Ish(G)$. We wish to create an Ish region $(w,\varepsilon)$ with ceiling partition $(\a,\b)$. To do this, we choose $w(1),\ldots,w(n)$ and then place $a_i$ circles above the symbol $b_i$. This will be an Ish ceiling diagram as long as the symbols $b_i$ occur in order, to the right of $1$. That is, the permutation $w$ must satisfy
\begin{equation}
\label{eq:abishcondition}
w^{-1}(1)<w^{-1}(b_1)<\cdots <w^{-1}(b_{n-k}).
\end{equation}
There are $\binom{n}{n-k+1}$ ways to place these symbols and then $(k-1)!$ ways to place the remaining symbols, proving (1). To prove (2), recall that $(w,\varepsilon)$ has $n-j+w^{-1}(1)$ degrees of freedom, where $j$ is the largest index such that $\varepsilon_j\neq 0$. In our case $j=w^{-1}(b_{n-k})$, so $w$ must satisfy the condition $w^{-1}(b_{n-k})-w^{-1}(1)=n-d$. First we can choose the pair $(w^{-1}(1),w^{-1}(b_{n-k}))$ in $d$ ways. Having done this, the rest of the permutation is subject  to \eqref{eq:abishcondition}. There are $\binom{n-d-1}{n-k-1}$ ways to place symbols $b_1,\ldots,b_{n-k-1}$ (left to right) in the $n-d-1$ positions between $1$ and $b_{n-k}$, and then there are $(k-1)!$ ways to place the remaining $k-1$ symbols. The result follows.
\end{proof}

Once again, note that this proof was more direct for Ish than for Shi. In fact, the calculation of formula \eqref{eq:candd} for Ish was the inspiration for Lemma \ref{lem:rhoades}.

\subsection{Concluding remarks}
The notion of a ceiling partition has independent interest, beyond the proof of Theorem \ref{thm:candd}. In particular, it leads to a new proof that the Shi arrangement $\Shi(n)$ has $(n+1)^{n-1}$ regions. Consider the collection of maps from $[n]$ into a set of size $x$. On one hand, there are $x^n$ of these. On the other hand, there are $\Stir(n,k)x(x-1)\cdots (x-k+1)$ such maps with image of size $k$, where $\Stir(n,k)$ is the number of partitions of $[n]$ into $k$ blocks (fibers). This proves the famous polynomial identity:
\begin{equation*}
x^n = \sum_{k=1}^n \Stir(n,k)x(x-1)\cdots (x-k+1).
\end{equation*}
Dividing by $x$ and substituting $x=n+1$ yields
\begin{equation*}
(n+1)^{n-1} = \sum_{k=1}^n \Stir(n,k) \frac{n!}{(n-k+1)!},
\end{equation*}
where the right hand side counts regions of $\Shi(n)$ by the number $k$ of blocks in their ceiling partition.

Finally, here are some problems for future research.

\begin{enumerate}
\item The original motivation for this paper was to find a bijection between regions of Shi and Ish. We solved this problem for dominant regions, but not in general. Based on Theorem \ref{thm:candd}, one should look for a bijection between Shi ceiling diagrams $(w,\pi)$ and Ish ceiling diagrams $(w,\varepsilon)$ with a fixed ceiling partition $(\a,\b)$ and $d$ degrees of freedom. Note that this bijection {\bf cannot} preserve the permutation $w$.
\medskip

\item Following Theorem \ref{th:charpoly}, find a direct bijection between points of the finite vector space $\FF_p^n$ in the complements of the Shi arrangement $\Shi(G)_p$ and the Ish arrangement $\Ish(G)_p$.
\medskip

\item The Shi arrangement is a famous example of a {\sf free hyperplane arrangement}. Investigate the freeness of Ish arrangements.
\medskip

\item Define and study an extended Ish arrangement corresponding to the {\sf extended Shi arrangement}:
\begin{equation*}\Shi(n,m):=\{ x_i-x_j = a : 1\leq i<j\leq n \,,\, -m+1\leq a\leq m\}.
\end{equation*}

\item To what extent do the results of this paper apply to other deformations of the Coxeter arrangement?
\medskip

\item The deleted Shi arrangements exist for arbitrary crystallographic reflection groups. Define and study Ish arrangements for other reflection groups. Ish arrangements were invented to study $q,t$-Catalan numbers; this feature should extend to other types.

\end{enumerate}

\section{Acknowledgements}
The authors are grateful to Christos Athanasiadis, Susanna Fishel, Christian Krattenthaler,
Vic Reiner, and Richard Stanley for helpful conversations.  This project arose out of
an AMS sectional meeting held at Penn State in October 2009.

\bibliography{../bib/my}

\end{document}